\documentclass[a4paper,12pt]{amsart}
\usepackage{graphicx,color}
\usepackage{pict2e}
\usepackage{amsmath}
\newtheorem{theorem}{Theorem}

\newtheorem{remark}{Remark}

\begin{document}
\title{Some chain maps on Khovanov complexes and Reidemeister moves}
\author{Noboru Ito}\thanks{The author is a Research Fellow of the Japan Society for the Promotion of Science.  This work was partly supported by KAKENHI}

\begin{abstract}
We introduce some chain maps between Khovanov complexes.  Each of the chain maps commutes with a chain homotopy map and a retraction maps which obtain a Reidemeister invariance of Khovanov homology.  
\end{abstract}

\maketitle

\section{Introduction.  }
Let $\mathcal{C}_{3}$ $=$ $\mathcal{C}\left(~
\begin{minipage}{60pt}
\begin{picture}(50,60)
\cbezier(0,0)(0,3)(7,7)(12,12)
\qbezier(19,19)(25,25)(31,31)
\cbezier(38,38)(45,45)(50,50)(50,56)
\cbezier(0,56)(0,50)(17.25,44.25)(34.5,34.5)
\cbezier(34.5,34.5)(44,30)(50,20)(50,0)
\cbezier(17,48)(21,53)(24,54)(24,56)
\cbezier(12,41)(-1,25)(24,10)(24,0)
\put(4,40){$c$}
\put(41,31){$b$}
\put(3,12){$a$}
\end{picture}
\end{minipage}\right)$ be a Khovanov complex of a link diagram and $\mathcal{C}_{2}$ $=$ $\mathcal{C}\left(~
\begin{minipage}{60pt}
\begin{picture}(50,60)
\cbezier(0,0)(0,3)(7,7)(12,12)
\qbezier(19,19)(25,25)(31,31)
\cbezier(38,38)(45,45)(50,50)(50,56)
\cbezier(14,43)(20,45)(17.25,44.25)(34.5,34.5)
\cbezier(34.5,34.5)(44,30)(50,20)(50,0)
\qbezier(24,56)(12,40)(0,56)
\cbezier(10,39)(1.5,25)(24,10)(24,0)
\qbezier(10,39)(12,42.3)(14,43)
\put(41,31){$b$}
\put(3,12){$a$}
\end{picture}
\end{minipage}
\right)$, $\mathcal{C}_{1}$ $=$ $\mathcal{C}\left(~
\begin{minipage}{60pt}
\begin{picture}(50,60)
\cbezier(0,0)(12,15)(20,20)(24,0)
\qbezier(19,19)(25,25)(31,31)
\cbezier(38,38)(45,45)(50,50)(50,56)
\cbezier(14,43)(20,45)(17.25,44.25)(34.5,34.5)
\cbezier(34.5,34.5)(44,30)(50,20)(50,0)
\qbezier(24,56)(12,40)(0,56)
\cbezier(10,39)(1.5,25)(15,15)(19,19)
\qbezier(10,39)(12,42.3)(14,43)
\put(41,31){$b$}
\end{picture}
\end{minipage}
\right)$, and $\mathcal{C}_{1'}$ $=$ $\mathcal{C}\left(~
\begin{minipage}{60pt}
\begin{picture}(50,60)
\cbezier(0,0)(0,3)(7,7)(12,12)
\qbezier(19,19)(25,25)(31,31)
\cbezier(38,36.5)(45,45)(50,50)(50,56)
\cbezier(38,31)(44,30)(50,20)(50,0)
\qbezier(38,31)(35,32.25)(38,36.5)
\cbezier(14,43)(20,45)(17.25,44.25)(31,34.5)
\qbezier(31,31)(32,32.25)(31,34.5)
\qbezier(24,56)(12,40)(0,56)
\cbezier(10,39)(1.5,25)(24,10)(24,0)
\qbezier(10,39)(12,42.3)(14,43)
\put(3,12){$a$}
\end{picture}
\end{minipage}
\right)$ be its subcomplexes.  
There exist chain homotopy maps $h_{J}$ ($J$ $=$ $1$, $1'$, $2$, $3$) relating the identity maps $\mathcal{C}_{J}$ $\to$ $\mathcal{C}_{J}$ with compositions $\operatorname{in} \circ \rho_{J}$ of an inclusion maps $\operatorname{in}$ and a retraction $\rho_{J}$ for Reidemeister moves $1$: $\begin{minipage}{60pt}
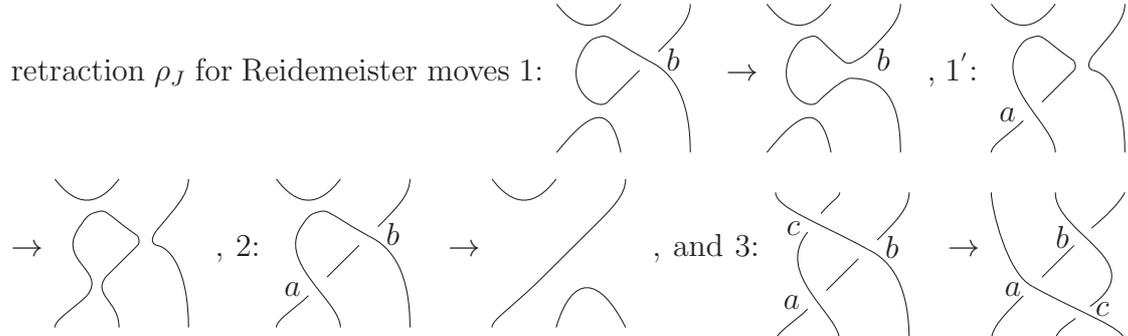
\begin{picture}(50,60)
\cbezier(0,0)(12,15)(20,20)(24,0)
\qbezier(19,19)(25,25)(31,31)
\cbezier(38,38)(45,45)(50,50)(50,56)
\cbezier(14,43)(20,45)(17.25,44.25)(34.5,34.5)
\cbezier(34.5,34.5)(44,30)(50,20)(50,0)
\qbezier(24,56)(12,40)(0,56)
\cbezier(10,39)(1.5,25)(15,15)(19,19)
\qbezier(10,39)(12,42.3)(14,43)
\put(41,31){$b$}
\end{picture}
\end{minipage}$ $\to$ $\begin{minipage}{60pt}
\begin{picture}(50,60)
\cbezier(0,0)(12,15)(20,20)(24,0)
\cbezier(19,19)(25,25)(29,27.5)(31,28)
\cbezier(38,38)(45,45)(50,50)(50,56)
\cbezier(38,38)(36,36)(34,35)(32,34)
\cbezier(14,43)(20,45)(17.25,44.25)(28,35)
\qbezier(28,35)(30,33)(32,34)
\cbezier(31,28)(44,28)(50,20)(50,0)
\qbezier(24,56)(12,40)(0,56)
\cbezier(10,39)(1.5,25)(15,15)(19,19)
\qbezier(10,39)(12,42.3)(14,43)
\put(41,31){$b$}
\end{picture}
\end{minipage}$, 
$1^{'}$: 
$\begin{minipage}{60pt}
\begin{picture}(50,60)
\cbezier(0,0)(0,3)(7,7)(12,12)
\qbezier(19,19)(25,25)(31,31)
\cbezier(38,36.5)(45,45)(50,50)(50,56)
\cbezier(38,31)(44,30)(50,20)(50,0)
\qbezier(38,31)(35,32.25)(38,36.5)
\cbezier(14,43)(20,45)(17.25,44.25)(31,34.5)
\qbezier(31,31)(32,32.25)(31,34.5)
\qbezier(24,56)(12,40)(0,56)
\cbezier(10,39)(1.5,25)(24,10)(24,0)
\qbezier(10,39)(12,42.3)(14,43)
\put(3,12){$a$}
\end{picture}
\end{minipage}$ $\to$ 
$\begin{minipage}{60pt}
\begin{picture}(50,60)
\cbezier(0,0)(0,3)(7,7)(12,12)
\qbezier(19,19)(25,25)(31,31)
\cbezier(38,36.5)(45,45)(50,50)(50,56)
\cbezier(38,31)(44,30)(50,20)(50,0)
\qbezier(38,31)(35,32.25)(38,36.5)
\cbezier(14,43)(20,45)(17.25,44.25)(31,34.5)
\qbezier(31,31)(32,32.25)(31,34.5)
\qbezier(24,56)(12,40)(0,56)
\cbezier(10,39)(-1,25)(20,22)(12,12)
\qbezier(10,39)(12,42.3)(14,43)
\cbezier(19,19)(14,12)(24,10)(24,0)
\end{picture}
\end{minipage}$, $2$: $\begin{minipage}{60pt}
\begin{picture}(50,60)
\cbezier(0,0)(0,3)(7,7)(12,12)
\qbezier(19,19)(25,25)(31,31)
\cbezier(38,38)(45,45)(50,50)(50,56)
\cbezier(14,43)(20,45)(17.25,44.25)(34.5,34.5)
\cbezier(34.5,34.5)(44,30)(50,20)(50,0)
\qbezier(24,56)(12,40)(0,56)
\cbezier(10,39)(-1,25)(24,10)(24,0)
\qbezier(10,39)(12,42.3)(14,43)
\put(41,31){$b$}
\put(3,12){$a$}
\end{picture}
\end{minipage}$ $\to$ 
$\begin{minipage}{60pt}
\begin{picture}(50,60)
\cbezier(0,0)(0,3)(7,7)(19,19)
\qbezier(19,19)(25,25)(38,38)
\cbezier(38,38)(45,45)(50,50)(50,56)
\cbezier(24,0)(30,20)(40,20)(50,0)
\qbezier(24,56)(12,40)(0,56)
\end{picture}
\end{minipage}$, and $3$: $\begin{minipage}{60pt}
\begin{picture}(50,70)
\cbezier(0,0)(0,3)(7,7)(12,12)
\qbezier(19,19)(25,25)(31,31)
\cbezier(38,38)(45,45)(50,50)(50,56)
\cbezier(0,56)(0,50)(17.25,44.25)(34.5,34.5)
\cbezier(34.5,34.5)(44,30)(50,20)(50,0)
\cbezier(17,48)(21,53)(24,54)(24,56)
\cbezier(12,41)(-1,25)(24,10)(24,0)
\put(4,40){$c$}
\put(41,31){$b$}
\put(3,12){$a$}
\end{picture}
\end{minipage}$ $\to$ $\begin{minipage}{60pt}
\begin{picture}(50,70)
\cbezier(37.5,42.5)(43,47)(48,52)(50,56)
\qbezier(31,36)(25,30)(19,24)
\cbezier(12,17)(5,10)(0,5)(0,0)
\cbezier(50,0)(50,5)(32.75,10.75)(15.5,20.5)
\cbezier(15.5,20.5)(6,25)(0,50)(0,56)
\cbezier(31.5,8.5)(30,7)(24,4)(24,0)
\cbezier(37,13)(61,30)(24,40)(24,56)
\put(39,9.5){$c$}
\put(5,17){$a$}
\put(24,35){$b$}
\end{picture}
\end{minipage}$ (Section \ref{right}, Appendix \ref{homotopy}, \ref{retraction}).  
This paper will be show that the natural chain maps $\pi_{J}$: $\mathcal{C}_{J}$ $\to$ $\mathcal{C}_{J-1}$ ($J$ $=$ $2$, $3$) satisfy the relations $h_{J-1} \circ \pi_{J}$ $=$ $\pi_{J} \circ h_{J}$ (Theorem \ref{commute}) and similar relations for $\rho_{J}$ (Theorem \ref{commute_rho}).  

In section 1 chain maps $\pi_{J}$ are defined.  In section 2 relations of $h_{J}$, $\rho_{J}$, and $\pi_{J}$ are given.  In section 3 a map $\tilde{\pi}_{2}$ similar to $\pi_{2}$ is introduced.  In section 4 we obtain the proof of the right twisted first Reidemeister invariance of Khovanov homology for a general differential.  Appendix contains the definitions $h_{J}$ and $\rho_{J}$ provided by \cite{ito3, ito5}.  All notations in this paper and the definition of the differential $\delta_{s, t}$ follows \cite{ito5}.  

\section{The chain map $\pi_{J}$.  }
$\pi_{3} : \mathcal{C}_{3} \to \mathcal{C}_{2}$ is defined by
\begin{equation}\label{pi3}
\begin{split}
\begin{minipage}{60pt}
\begin{picture}(50,70)
            \cbezier(0,0)(0,3)(7,7)(12,12)
            \qbezier(19,19)(25,25)(31,31)
            \cbezier(38,38)(45,45)(50,50)(50,56)
            \cbezier(0,56)(0,50)(17.25,44.25)(34.5,34.5)
            \cbezier(34.5,34.5)(44,30)(50,20)(50,0)
            \cbezier(17,48)(21,53)(24,54)(24,56)
            \cbezier(12,41)(-1,25)(24,10)(24,0)
            \put(0,40){$c$}
            \put(41,31){$b$}
            \put(3,12){$a$}
            {\color{blue}{\put(8.5,44.5){\circle*{3}}
            \put(20.5,44.5){\circle*{3}}
            \put(9.5,44.5){\line(1,0){10}}
            }}
        \end{picture}
    \end{minipage} &\mapsto \begin{minipage}{60pt}
\begin{picture}(50,60)
\cbezier(0,0)(0,3)(7,7)(12,12)
\qbezier(19,19)(25,25)(31,31)
\cbezier(38,38)(45,45)(50,50)(50,56)
\cbezier(14,43)(20,45)(17.25,44.25)(34.5,34.5)
\cbezier(34.5,34.5)(44,30)(50,20)(50,0)
\qbezier(24,56)(12,40)(0,56)
\cbezier(10,39)(1.5,25)(24,10)(24,0)
\qbezier(10,39)(12,42.3)(14,43)
\put(41,31){$b$}
\put(3,12){$a$}
\end{picture}
\end{minipage}, \\
    \begin{minipage}{60pt}
        \begin{picture}(50,70)
            \cbezier(0,0)(0,3)(7,7)(12,12)
            \qbezier(19,19)(25,25)(31,31)
            \cbezier(38,38)(45,45)(50,50)(50,56)
            \cbezier(0,56)(0,50)(17.25,44.25)(34.5,34.5)
            \cbezier(34.5,34.5)(44,30)(50,20)(50,0)
            \cbezier(17,48)(21,53)(24,54)(24,56)
            \cbezier(12,41)(-1,25)(24,10)(24,0)
            {\color{red}{\put(14,49){\circle*{3}}
\put(14,39){\circle*{3}}
\put(14,39){\line(0,1){10}}}}
            \put(4,40){$c$}
            \put(41,31){$b$}
            \put(3,12){$a$}
        \end{picture}
    \end{minipage}
    &\mapsto 0.  
\end{split}
\end{equation}

$\pi_{2} : \mathcal{C}_{2} \to \mathcal{C}_{1}$ is defined by 
\begin{equation}\label{pi2}
\begin{split}
\begin{minipage}{60pt}
\begin{picture}(50,60)
\cbezier(0,0)(0,3)(7,7)(12,12)
\qbezier(19,19)(25,25)(31,31)
\cbezier(38,38)(45,45)(50,50)(50,56)
\cbezier(14,43)(20,45)(17.25,44.25)(34.5,34.5)
\cbezier(34.5,34.5)(44,30)(50,20)(50,0)
\qbezier(24,56)(12,40)(0,56)
\cbezier(10,39)(1.5,25)(24,10)(24,0)
\qbezier(10,39)(12,42.3)(14,43)
\put(41,31){$b$}
\put(3,12){$a$}
{\color{blue}{\put(11,15){\circle*{3}}
\put(21,15){\circle*{3}}
\put(11,15){\line(1,0){10}}}}
\end{picture}
\end{minipage}
    &\mapsto \begin{minipage}{60pt}
\begin{picture}(50,60)
\cbezier(0,0)(12,15)(20,20)(24,0)
\qbezier(19,19)(25,25)(31,31)
\cbezier(38,38)(45,45)(50,50)(50,56)
\cbezier(14,43)(20,45)(17.25,44.25)(34.5,34.5)
\cbezier(34.5,34.5)(44,30)(50,20)(50,0)
\qbezier(24,56)(12,40)(0,56)
\cbezier(10,39)(1.5,25)(15,15)(19,19)
\qbezier(10,39)(12,42.3)(14,43)
\put(41,31){$b$}
\end{picture}
\end{minipage}  \\
\begin{minipage}{60pt}
\begin{picture}(50,70)
\cbezier(0,0)(0,3)(7,7)(12,12)
\qbezier(19,19)(25,25)(31,31)
\cbezier(38,38)(45,45)(50,50)(50,56)
\cbezier(0,56)(0,50)(17.25,44.25)(34.5,34.5)
\cbezier(34.5,34.5)(44,30)(50,20)(50,0)
\cbezier(17,48)(21,53)(24,54)(24,56)
\cbezier(12,41)(-1,25)(24,10)(24,0)
\put(4,40){$c$}
\put(41,31){$b$}
\put(3,12){$a$}
{\color{red}{\put(15,21){\circle*{3}}
\put(15,10){\circle*{3}}
\put(15,11){\line(0,1){10}}}}
\end{picture}
\end{minipage}
    &\mapsto 0.  
\end{split}
\end{equation}    
\begin{theorem}\label{commute}
The chain maps $\pi_{3}$ and $\pi_{2}$ satisfy the following.  
\begin{align}
h_{2} \circ \pi_{3} &= \pi_{3} \circ h_{3},  \label{commute_1}\\
h_{1} \circ \pi_{2} &= \pi_{2} \circ h_{2}.  \label{commute_2}
\end{align}
\end{theorem}
\begin{proof}
$\delta_{s, t} \circ \pi_{J}$ $=$ $\pi_{J} \circ \delta_{s, t}$ and $h_{J-1} \circ \pi_{J}$ $=$ $\pi_{J} \circ h_{J}$  ($J$ $=$ $2$, $3$) are proved by direct computation for every generator of $\mathcal{C}_{J}$.  
\end{proof}

Let $\mathcal{C}'_{3} = \mathcal{C}\left(~\begin{minipage}{60pt}
\begin{picture}(50,70)
\cbezier(0,0)(0,3)(7,7)(12,12)
\qbezier(19,19)(25,25)(31,31)
\cbezier(38,38)(45,45)(50,50)(50,56)
\cbezier(0,56)(0,50)(17.25,44.25)(34.5,34.5)
\cbezier(34.5,34.5)(44,30)(50,20)(50,0)
\cbezier(17,48)(21,53)(24,54)(24,56)
\cbezier(12,41)(-1,25)(24,10)(24,0)
            {\color{red}{\put(16,21){\circle*{3}}
\put(16,10){\circle*{3}}
\put(16,11){\line(0,1){10}}}}
{\color{blue}
            {\put(4.5,44.5){\circle*{3}}
\put(16.5,44.5){\circle*{3}}
\put(5.5,44.5){\line(1,0){10}}
            }}
            {\color{blue}
            {\put(32.5,34.5){\circle*{3}}
\put(20.5,34.5){\circle*{3}}
\put(21.5,34.5){\line(1,0){10}}
            }}
            \put(13,14){$q$}
            \put(45,45){$p$}
            \put(4,48){$r$}
        \end{picture}
    \end{minipage} \otimes [xa] + \begin{minipage}{60pt}
\begin{picture}(50,70)
\cbezier(0,0)(0,3)(7,7)(12,12)
\qbezier(19,19)(25,25)(31,31)
\cbezier(38,38)(45,45)(50,50)(50,56)
\cbezier(0,56)(0,50)(17.25,44.25)(34.5,34.5)
\cbezier(34.5,34.5)(44,30)(50,20)(50,0)
\cbezier(17,48)(21,53)(24,54)(24,56)
\cbezier(12,41)(-1,25)(24,10)(24,0)
            \put(15,28){\text{$-$}}
            {\color{blue}
            {\put(8,44.5){\circle*{3}}
\put(20,44.5){\circle*{3}}
\put(9,44.5){\line(1,0){10}}
            }}
            {\color{blue}{\put(8,15){\circle*{3}}
\put(18,15){\circle*{3}}
\put(9,15){\line(1,0){10}}}}
{\color{red}{\put(27,39){\circle*{3}}
\put(27,29){\circle*{3}}
\put(27,29){\line(0,1){10}}}}
\put(-6,-5){$p:q$}
\put(45,0){$q:p$}
\put(5,50){$\tilde{r}$}
        \end{picture}
    \end{minipage}\!\!\!\! \otimes [xb], \begin{minipage}{60pt}
\begin{picture}(50,70)
\cbezier(0,0)(0,3)(7,7)(12,12)
\qbezier(19,19)(25,25)(31,31)
\cbezier(38,38)(45,45)(50,50)(50,56)
\cbezier(0,56)(0,50)(17.25,44.25)(34.5,34.5)
\cbezier(34.5,34.5)(44,30)(50,20)(50,0)
\cbezier(17,48)(21,53)(24,54)(24,56)
\cbezier(12,41)(-1,25)(24,10)(24,0)
            {\color{red}{\put(14,49){\circle*{3}}
\put(14,39){\circle*{3}}
\put(14,39){\line(0,1){10}}}}
        \end{picture}
    \end{minipage} \otimes [x]\right)$ and
    \[\]
     $\mathcal{C}'_{2} = \mathcal{C}\left(~\begin{minipage}{40pt}
        \begin{picture}(30,40)
           \qbezier(0,5)(40,20)(0,35)
            \qbezier(11,26)(1,19)(11,14)
            \qbezier(18,29)(24,32)(30,35)
            \qbezier(18,10)(24,7.5)(30,5)
  {\color{red}{\put(14,33){\circle*{3}}
\put(14,24){\circle*{3}}
\put(14,24){\line(0,1){8}}}}
            {\color{blue}{\put(10,17){\circle*{3}}
\put(10,8){\circle*{3}}
\put(10,9){\line(0,1){8}}}}
\put(-5,19){$p$}
\put(18,19){$q$}
        \end{picture}
    \end{minipage} \otimes [xa] + \begin{minipage}{40pt}
        \begin{picture}(30,40)
            \qbezier(0,5)(40,20)(0,35)
            \qbezier(11,26)(1,19)(11,14)
            \qbezier(18,29)(24,32)(30,35)
            \qbezier(18,10)(24,7.5)(30,5)
            {\color{blue}
            {\put(8,27.5){\circle*{3}}
\put(20,27.5){\circle*{3}}
\put(9,27.5){\line(1,0){10}}
            }}
            {\color{red}{\put(16,12){\circle*{3}}
\put(5,12){\circle*{3}}
\put(5,12){\line(1,0){10}}}}
\put(5,17){\text{$-$}}
\put(1,37){$p:q$}
\put(0,1){$q:p$}
        \end{picture}
    \end{minipage} \otimes [xb]~\right)$.  $\mathcal{C}'_{J}$ is a subcomplex of $\mathcal{C}_{J}$ ($J$ $=$ $2$, $3$).  
We define $\pi'_{3} : \mathcal{C}'_{3} \to \mathcal{C}_{2}$  by (\ref{pi3}).  

\begin{theorem}\label{commute_rho}
\begin{align}
\rho_{2} \circ \pi_{3} &= \pi'_{3} \circ \rho_{3}.  
\end{align}
\end{theorem}
\begin{proof}
$\delta_{s, t} \circ \pi'_{3}$ $=$ $\pi'_{3} \circ \delta_{s, t}$ and $\rho_{2} \circ \pi_{3}$ $=$ $\pi'_{3} \circ \rho_{3}$ are proved by direct computation for every generator of $\mathcal{C}_{3}$.  
\end{proof}
\section{A similar map $\tilde{\pi}_{2}$ to $\pi_{2}$.  }
In this section we introduce a map $\tilde{\pi}_{2}$.  It is not chain maps, but it has similar property $h_{1'} \circ \tilde{\pi}_{2}$ $=$ $\tilde{\pi}_{2} \circ h_{1'}$ (Theorem \ref{similar}).  

The map $\tilde{\pi}_{2} : \mathcal{C}_{2} \to \mathcal{C}_{1'}$ is defined by 
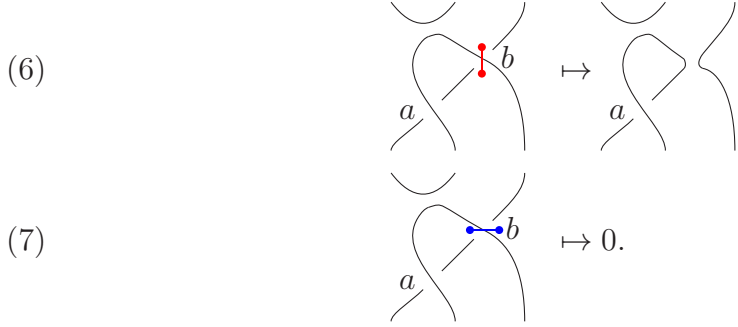
\begin{align}
\begin{minipage}{60pt}
\begin{picture}(50,60)
\cbezier(0,0)(0,3)(7,7)(12,12)
\qbezier(19,19)(25,25)(31,31)
\cbezier(38,38)(45,45)(50,50)(50,56)
\cbezier(14,43)(20,45)(17.25,44.25)(34.5,34.5)
\cbezier(34.5,34.5)(44,30)(50,20)(50,0)
\qbezier(24,56)(12,40)(0,56)
\cbezier(10,39)(1.5,25)(24,10)(24,0)
\qbezier(10,39)(12,42.3)(14,43)
\put(41,31){$b$}
\put(3,12){$a$}
{\color{red}{\put(34,39){\circle*{3}}
\put(34,29){\circle*{3}}
\put(34,29){\line(0,1){10}}
}}
\end{picture}
\end{minipage} &\mapsto \begin{minipage}{60pt}
\begin{picture}(50,60)
\cbezier(0,0)(0,3)(7,7)(12,12)
\qbezier(19,19)(25,25)(31,31)
\cbezier(38,36.5)(45,45)(50,50)(50,56)
\cbezier(38,31)(44,30)(50,20)(50,0)
\qbezier(38,31)(35,32.25)(38,36.5)
\cbezier(14,43)(20,45)(17.25,44.25)(31,34.5)
\qbezier(31,31)(32,32.25)(31,34.5)
\qbezier(24,56)(12,40)(0,56)
\cbezier(10,39)(1.5,25)(24,10)(24,0)
\qbezier(10,39)(12,42.3)(14,43)
\put(3,12){$a$}
\end{picture}
\end{minipage}\\
\begin{minipage}{60pt}
\begin{picture}(50,60)
\cbezier(0,0)(0,3)(7,7)(12,12)
\qbezier(19,19)(25,25)(31,31)
\cbezier(38,38)(45,45)(50,50)(50,56)
\cbezier(14,43)(20,45)(17.25,44.25)(34.5,34.5)
\cbezier(34.5,34.5)(44,30)(50,20)(50,0)
\qbezier(24,56)(12,40)(0,56)
\cbezier(10,39)(1.5,25)(24,10)(24,0)
\qbezier(10,39)(12,42.3)(14,43)
\put(43,31){$b$}
\put(3,12){$a$}
{\color{blue}{\put(40.5,34.5){\circle*{3}}
\put(29.5,34.5){\circle*{3}}
\put(30.5,34.5){\line(1,0){10}}
}}
\end{picture}
\end{minipage} &\mapsto 0.  
\end{align}
\begin{theorem}\label{similar}
\begin{equation}
h_{1'} \circ \tilde{\pi}_{2} = \tilde{\pi}_{2} \circ h_{1'}.  
\end{equation}
\end{theorem}
\begin{proof}
$h_{1'} \circ \tilde{\pi}_{2}$ $=$ $\tilde{\pi}_{2} \circ h_{1'}$ is proved by direct computation for every generator of $\mathcal{C}_{2}$.  
\end{proof}

\section{Right twisted first Reidemeister invariance for the general differential $\delta_{s, t}$.  }\label{right}
In this section we will show that the right twisted first Reidemeister invariance of Khovanov homology  because the proof of this case is missing in \cite{ito5}.  

The right twisted first Reidemeister move is $D'$ $=$ $~\begin{minipage}{30pt}
        \begin{picture}(30,30)
        \put(-3,12){$a$}
\qbezier(6.6,12)(3.3,8.5)(0,5)
\qbezier(0,25)(20,-5)(25,12)
\qbezier(11,16)(20,25)(24.5,17)
\qbezier(24.5,17)(25.5,14.5)(25,12)
        \end{picture}
    \end{minipage}$ $\stackrel{1}{\sim}$ $\begin{minipage}{30pt}
        \begin{picture}(30,30)
        \cbezier(0,3)(30,5)(30,25)(0,27) 
        \end{picture}
    \end{minipage} $ $=$ $D$, we consider the composition 

\begin{equation}
\mathcal{C}(D') = \mathcal{C} \oplus \mathcal{C}_{\rm{contr}} \stackrel{\rho_{1}}{\to} \mathcal{C} \stackrel{\operatorname{isom}}{\to} \mathcal{C}(D)
\end{equation} 
where $a$ is a crossing and $\mathcal{C}$, $\mathcal{C}_{\rm{contr}}$, $\rho_{1}$ and the isomorphism are defined in the following formulas (\ref{1st-sum})--(\ref{1st-iso}).  

First, 
\begin{equation}\label{1st-sum}
\begin{split}
&\mathcal{C} := \mathcal{C} \left(~
    \begin{minipage}{30pt}
        \begin{picture}(30,30)
\qbezier(6.6,12)(3.3,8.5)(0,5)
\qbezier(0,25)(20,-5)(25,12)
\qbezier(11,16)(20,25)(24.5,17)
\qbezier(24.5,17)(25.5,14.5)(25,12)
{\color{red}{\put(8,18){\circle*{3}}
\put(8,9){\circle*{3}}
\put(8,10){\line(0,1){8}}}}
\put(-3.5,12.5){\text{$p$}}
\put(14,11.5){\text{$-$}}
        \end{picture}
    \end{minipage}\!\! \otimes [x] 
   \right), \\
   &\mathcal{C}_{\rm{contr}} := \mathcal{C}\left(~\qquad 
\begin{minipage}{30pt}
        \begin{picture}(30,30)
\qbezier(6.6,13)(3.3,10)(0,6)
\qbezier(0,26)(20,-4)(25,13)
\qbezier(11,17)(20,26)(24.5,18)
\qbezier(24.5,18)(25.5,16)(25,13)
{\color{red}{\put(8,19.5){\circle*{3}}
\put(8,10.5){\circle*{3}}
\put(8,11.5){\line(0,1){8}}}}
\put(-20.5,12.5){\text{$p:p$}}
\put(25,25){\text{$p:p$}}
\qbezier(18,15)(15,30)(24.5,27)
        \end{picture}
    \end{minipage} \!\! \otimes [x],   
    \begin{minipage}{30pt}
        \begin{picture}(30,30)
\qbezier(6.6,13.5)(3.3,10.5)(0,6.5)
\qbezier(0,26.5)(20,-3.5)(25,13.5)
\qbezier(11,17.5)(20,26.5)(24.5,18.5)
\qbezier(24.5,18.5)(25.5,16.5)(25,13.5)
{\color{blue}{\put(12.5,15.5){\circle*{3}}
\put(4.5,15.5){\circle*{3}}
\put(4.5,15.5){\line(1,0){8}}}}
\put(15.5,14){\text{$p$}}
        \end{picture}
    \end{minipage} \!\! \otimes [xa] \right).  
\end{split}
\end{equation}

Second, the retraction 
$\rho_{1} : \mathcal{C}\left(~
    \begin{minipage}{30pt}
        \begin{picture}(30,30)
\qbezier(6.6,12)(3.3,8.5)(0,5)
\qbezier(0,25)(20,-5)(25,12)
\qbezier(11,16)(20,25)(24.5,17)
\qbezier(24.5,17)(25.5,14.5)(25,12)
        \end{picture}
    \end{minipage}
\right)$ $\to$ 
$ 
\mathcal{C}\left(~\begin{minipage}{30pt}
        \begin{picture}(30,30)
\qbezier(6.6,13)(3.3,9.5)(0,6)
\qbezier(0,26)(20,-4)(25,13)
\qbezier(11,17)(20,26)(24.5,18)
\qbezier(24.5,18)(25.5,15.5)(25,13)
{\color{red}{\put(8,19.5){\circle*{3}}
\put(8,10.5){\circle*{3}}
\put(8,11.5){\line(0,1){8}}}}
\put(-3.5,12.5){\text{$p$}}
\put(14,12.5){\text{$-$}}
        \end{picture}
    \end{minipage}\!\! \otimes [x] \right)$
is defined by the formulas

\begin{equation}\label{1st-ret}
\begin{split}
&\begin{minipage}{30pt}
        \begin{picture}(30,30)
        \qbezier(6.6,13)(3.3,9.5)(0,6)
        \qbezier(0,26)(20,-4)(25,13)
        \qbezier(11,17)(20,26)(24.5,18)
        \qbezier(24.5,18)(25.5,15.5)(25,12)
{\color{red}{\put(8,19.5){\circle*{3}}
\put(8,10.5){\circle*{3}}
\put(8,11.5){\line(0,1){8}}}}
\put(-3.5,12.5){\text{$p$}}
\put(14,12.5){\text{$-$}}
        \end{picture}
    \end{minipage} \!\! \otimes [x] \mapsto \ 
\begin{minipage}{30pt}
        \begin{picture}(30,30)
\qbezier(6.6,13)(3.3,9.5)(0,6)
        \qbezier(0,26)(20,-4)(25,13)
        \qbezier(11,17)(20,26)(24.5,18)
        \qbezier(24.5,18)(25.5,15.5)(25,12)
{\color{red}{\put(8,19.5){\circle*{3}}
\put(8,10.5){\circle*{3}}
\put(8,11.5){\line(0,1){8}}}}
\put(-3.5,12.5){\text{$p$}}
\put(14,12.5){\text{$-$}}
        \end{picture}
    \end{minipage}\!\! \otimes [x], \\
&    \begin{minipage}{30pt}
        \begin{picture}(30,30)
\qbezier(6.6,13)(3.3,9.5)(0,6)
        \qbezier(0,26)(20,-4)(25,13)
        \qbezier(11,17)(20,26)(24.5,18)
        \qbezier(24.5,18)(25.5,15.5)(25,12)
{\color{red}{\put(8,19.5){\circle*{3}}
\put(8,10.5){\circle*{3}}
\put(8,11.5){\line(0,1){8}}}}
\put(-3.5,12.5){\text{$p$}}
\put(14,12.5){\text{$+$}}
        \end{picture}
    \end{minipage} \!\! \otimes [x] \mapsto \  \begin{minipage}{30pt}
\begin{picture}(30,30)
\qbezier(6.6,13)(3.3,10)(0,6)
\qbezier(0,26)(20,-4)(25,13)
\qbezier(24.5,18)(25.5,16)(25,13)
\qbezier(11,17)(20,26)(24.5,18)
{\color{red}{\put(8,19.5){\circle*{3}}
\put(8,10.5){\circle*{3}}
\put(8,11.5){\line(0,1){8}}}}
\put(-3.5,12.5){\text{$p$}}
\put(14,11.5){\text{$+$}}
\end{picture}
\end{minipage}\!\! \otimes [x] - \qquad \begin{minipage}{30pt}
\begin{picture}(30,30)
\qbezier(6.6,13)(3.3,10)(0,6)
\qbezier(0,26)(20,-4)(25,13)
\qbezier(24.5,18)(25.5,16)(25,13)
\qbezier(11,17)(20,26)(24.5,18)
{\color{red}{\put(8,19.5){\circle*{3}}
\put(8,10.5){\circle*{3}}
\put(8,11.5){\line(0,1){8}}}}
\qbezier(18,15)(15,30)(24.5,27)
\put(-20.5,12.5){\text{$p:p$}}
\put(25,25){\text{$p:p$}}
\end{picture}
\end{minipage}\!\! \otimes [x],\\   
    &\begin{minipage}{30pt}
        \begin{picture}(30,30)
\qbezier(6.6,13)(3.3,9.5)(0,6)
        \qbezier(0,26)(20,-4)(25,13)
        \qbezier(11,17)(20,26)(24.5,18)
        \qbezier(24.5,18)(25.5,15.5)(25,12)
{\color{blue}{\put(12.5,15.5){\circle*{3}}
\put(4.5,15.5){\circle*{3}}
\put(4.5,15.5){\line(1,0){8}}}}
\put(16,14){\text{$p$}}
        \end{picture}
    \end{minipage} \!\! \otimes [xa] \mapsto \  0.  
\end{split}
\end{equation}
We can verify that $\delta_{s, t} \circ \rho_{1}$ $=$ $\rho_{1} \circ \delta_{s, t}$.  Then $\rho_{1}$ is a chain map.  

Third, the isomorphism 
\begin{equation}\label{1st-iso1}
\begin{split}
&\mathcal{C}\left(~
    \begin{minipage}{30pt}
        \begin{picture}(30,30)
\qbezier(6.6,13)(3.3,9.5)(0,6)
        \qbezier(0,26)(20,-4)(25,13)
        \qbezier(11,17)(20,26)(24.5,18)
        \qbezier(24.5,18)(25.5,15.5)(25,12)
{\color{red}{\put(8,19.5){\circle*{3}}
\put(8,10.5){\circle*{3}}
\put(8,11.5){\line(0,1){8}}}}
\put(-3.5,12.5){\text{$p$}}
\put(14,12.5){\text{$-$}}
        \end{picture}
    \end{minipage}\!\! \otimes [x]  \right) 
    \to 
    \mathcal{C}\left(~
    \begin{minipage}{30pt}
        \begin{picture}(30,30)
        \cbezier(0,3)(30,5)(30,25)(0,27) 
        \end{picture}
    \end{minipage} \!\! \otimes [x]
    \right) 
    \end{split}
        \end{equation}
    is defined by the formulas 
\begin{equation}\label{1st-iso}
\begin{split}
&\begin{minipage}{30pt}
        \begin{picture}(30,30)
\qbezier(6.6,13)(3.3,9.5)(0,6)
        \qbezier(0,26)(20,-4)(25,13)
        \qbezier(11,17)(20,26)(24.5,18)
        \qbezier(24.5,18)(25.5,15.5)(25,12)
        {\color{red}{\put(8,19.5){\circle*{3}}
\put(8,10.5){\circle*{3}}
\put(8,11.5){\line(0,1){8}}}}
\put(-3.5,12.5){\text{$p$}}
\put(14,12.5){\text{$-$}}
        \end{picture}
    \end{minipage}\!\! \otimes [x] \mapsto~\ 
    \begin{minipage}{30pt}
        \begin{picture}(30,30)
        \cbezier(0,3)(30,5)(30,25)(0,27) 
        \put(5,12.5){\text{$p$}}
        \end{picture}
    \end{minipage} \!\! \otimes [x].  
        \end{split}
        \end{equation}  
The homotopy connecting $\operatorname{in}$ $\circ$ $\rho_{1}$ to the identity $:$ $\mathcal{C}\left(~
    \begin{minipage}{30pt}
        \begin{picture}(30,30)
\qbezier(6.6,13)(3.3,9.5)(0,6)
        \qbezier(0,26)(20,-4)(25,13)
        \qbezier(11,17)(20,26)(24.5,18)
        \qbezier(24.5,18)(25.5,15.5)(25,12)
                \end{picture}
    \end{minipage}
\right)$ $\to$ 
$\mathcal{C}\left(~
    \begin{minipage}{30pt}
        \begin{picture}(30,30)
\qbezier(6.6,13)(3.3,9.5)(0,6)
        \qbezier(0,26)(20,-4)(25,13)
        \qbezier(11,17)(20,26)(24.5,18)
        \qbezier(24.5,18)(25.5,15.5)(25,12)
                \end{picture}
    \end{minipage}
\right)$ such that $\delta_{s, t}$ $\circ$ $h_{1}$ $+$  $h_{1}$ $\circ$ $\delta_{s, t}$ $=$ $\operatorname{id} - \operatorname{in} \circ \rho_{1}$, is defined by the formulas: 
\begin{equation}\label{first-right-hom}
\begin{split}
\begin{minipage}{30pt}
        \begin{picture}(30,30)
\qbezier(6.6,13)(3.3,9.5)(0,6)
        \qbezier(0,26)(20,-4)(25,13)
        \qbezier(11,17)(20,26)(24.5,18)
        \qbezier(24.5,18)(25.5,15.5)(25,12)
{\color{red}{\put(8,19.5){\circle*{3}}
\put(8,10.5){\circle*{3}}
\put(8,11.5){\line(0,1){8}}}}
\put(-3.5,12.5){\text{$p$}}
\put(14,12.5){\text{$-$}}
        \end{picture}
    \end{minipage} \!\! \otimes [xa] &\mapsto~ 
    \begin{minipage}{30pt}
        \begin{picture}(30,30)
\qbezier(6.6,13)(3.3,9.5)(0,6)
        \qbezier(0,26)(20,-4)(25,13)
        \qbezier(11,17)(20,26)(24.5,18)
        \qbezier(24.5,18)(25.5,15.5)(25,12)
{\color{blue}{\put(12.5,15){\circle*{3}}
\put(4.5,15){\circle*{3}}
\put(4.5,15){\line(1,0){8}}}}
\put(15.5,12.5){\text{$p$}}
        \end{picture}
    \end{minipage} \!\! \otimes [x],~
    {\text{otherwise}} \mapsto 0.  
\end{split}
\end{equation}
\begin{remark}
The explicit formula (\ref{first-right-hom}) of the homotopy map $h_{1}$ in the case ($s$ $=$ $t$ $=$ $0$) of the original Khovanov homology is given by Oleg Viro \cite[Subsection 5.5]{viro}.  
\end{remark}
We can verify $\delta_{s, t}$ $\circ$ $h_{1}$ $+$  $h_{1}$ $\circ$ $\delta_{s, t}$ $=$ $\operatorname{id} - \operatorname{in} \circ \rho_{1}$ by a direct computation as follows.  
\begin{equation}\label{first-begin}
\begin{split}
\left(h_{1} \circ \delta_{s, t} + \delta_{s, t} \circ h_{1} \right)\left(
\begin{minipage}{30pt}
        \begin{picture}(30,30)
\qbezier(6.6,13)(3.3,9.5)(0,6)
        \qbezier(0,26)(20,-4)(25,13)
        \qbezier(11,17)(20,26)(24.5,18)
        \qbezier(24.5,18)(25.5,15.5)(25,12)
{\color{blue}{\put(12.5,15){\circle*{3}}
        \put(4.5,15){\circle*{3}}
        \put(4.5,15){\line(1,0){8}}}}
        \put(17.5,12.5){\text{$p$}}
        \end{picture}
    \end{minipage} \!\! \otimes [x]\right) 
    &= h_{1}\left(\!\!\!\qquad\begin{minipage}{30pt}
\begin{picture}(30,30)
\qbezier(6.6,13)(3.3,10)(0,6)
\qbezier(0,26)(20,-4)(25,13)
\qbezier(24.5,18)(25.5,16)(25,13)
\qbezier(11,17)(20,26)(24.5,18)
{\color{red}{\put(8,19.5){\circle*{3}}
\put(8,10.5){\circle*{3}}
\put(8,11.5){\line(0,1){8}}}}
\qbezier(18,15)(15,30)(24.5,27)
\put(-20.5,12.5){\text{$p:p$}}
\put(25,25){\text{$p:p$}}
\end{picture}
\end{minipage}\!\! \otimes [xa]
      \right)\\
    &= \begin{minipage}{30pt}
    \begin{picture}(30,30)
    \qbezier(6.6,13)(3.3,9.5)(0,6)
    \qbezier(0,26)(20,-4)(25,13)
    \qbezier(11,17)(20,26)(24.5,18)
    \qbezier(24.5,18)(25.5,15.5)(25,12)
    {\color{blue}{\put(13,15){\circle*{3}}
    \put(5,15){\circle*{3}}
    \put(5,15){\line(1,0){8}}}}
    \put(17.5,12.5){\text{$p$}}
    \end{picture}
    \end{minipage}\!\! \otimes [x]
    \\
    &= \left( \operatorname{id} - \rho_{1} \right)\left(
    \begin{minipage}{30pt}
        \begin{picture}(30,30)
\qbezier(6.6,13)(3.3,9.5)(0,6)
        \qbezier(0,26)(20,-4)(25,13)
        \qbezier(11,17)(20,26)(24.5,18)
        \qbezier(24.5,18)(25.5,15.5)(25,12)
{\color{blue}{\put(12.5,15){\circle*{3}}
        \put(4.5,15){\circle*{3}}
        \put(4.5,15){\line(1,0){8}}}}
        \put(17.5,12.5){\text{$p$}}
        \end{picture}
    \end{minipage} \!\! \otimes [x]
    \right).  
    \end{split}
    \end{equation}
    Similarly, 
  \begin{equation}
\begin{split}
\left(h_{1} \circ \delta_{s, t} + \delta_{s, t} \circ h_{1} \right)\left(
\begin{minipage}{30pt}
        \begin{picture}(30,30)
\qbezier(6.6,13)(3.3,9.5)(0,6)
        \qbezier(0,26)(20,-4)(25,13)
        \qbezier(11,17)(20,26)(24.5,18)
        \qbezier(24.5,18)(25.5,15.5)(25,12)
{\color{red}{\put(8,19.5){\circle*{3}}
\put(8,10.5){\circle*{3}}
\put(8,11.5){\line(0,1){8}}}}
\put(-3.5,12.5){\text{$p$}}
\put(14,12.5){\text{$+$}}
        \end{picture}
    \end{minipage} \!\! \otimes [x]\right) 
    &= \qquad \begin{minipage}{30pt}
\begin{picture}(30,30)
\qbezier(6.6,13)(3.3,10)(0,6)
\qbezier(0,26)(20,-4)(25,13)
\qbezier(24.5,18)(25.5,16)(25,13)
\qbezier(11,17)(20,26)(24.5,18)
{\color{red}{\put(8,19.5){\circle*{3}}
\put(8,10.5){\circle*{3}}
\put(8,11.5){\line(0,1){8}}}}
\qbezier(18,15)(15,30)(24.5,27)
\put(-20.5,12.5){\text{$p:p$}}
\put(25,25){\text{$p:p$}}
\end{picture}
\end{minipage}\!\! \otimes [x]\\
    &= \left( \operatorname{id} - \rho_{1} \right)\left(
    \begin{minipage}{30pt}
        \begin{picture}(30,30)
\qbezier(6.6,13)(3.3,9.5)(0,6)
        \qbezier(0,26)(20,-4)(25,13)
        \qbezier(11,17)(20,26)(24.5,18)
        \qbezier(24.5,18)(25.5,15.5)(25,12)
{\color{red}{\put(8,19.5){\circle*{3}}
\put(8,10.5){\circle*{3}}
\put(8,11.5){\line(0,1){8}}}}
\put(-3.5,12.5){\text{$p$}}
\put(14,12.5){\text{$+$}}
        \end{picture}
    \end{minipage} \!\! \otimes [x]
    \right),   
    \end{split}
    \end{equation}
        \begin{equation}\label{1st-last}
\left(h_{1} \circ \delta_{s, t} + \delta_{s, t} \circ h_{1} \right)\left(~\begin{minipage}{30pt}
        \begin{picture}(30,30)
\qbezier(6.6,13)(3.3,9.5)(0,6)
        \qbezier(0,26)(20,-4)(25,13)
        \qbezier(11,17)(20,26)(24.5,18)
        \qbezier(24.5,18)(25.5,15.5)(25,12)
        {\color{red}{\put(8,19.5){\circle*{3}}
\put(8,10.5){\circle*{3}}
\put(8,11.5){\line(0,1){8}}}}
\put(-5,14){\text{$p$}}
\put(14,12){\text{$-$}}
        \end{picture}
    \end{minipage} \!\! \otimes [x]
\right) 
    = 0
    = \left(\operatorname{id} - \rho_{1} \right)\left(~
    \begin{minipage}{30pt}
        \begin{picture}(30,30)
\qbezier(6.6,13)(3.3,9.5)(0,6)
        \qbezier(0,26)(20,-4)(25,13)
        \qbezier(11,17)(20,26)(24.5,18)
        \qbezier(24.5,18)(25.5,15.5)(25,12)
{\color{red}{\put(8,19.5){\circle*{3}}
\put(8,10.5){\circle*{3}}
\put(8,11.5){\line(0,1){8}}}}
\put(-5,14){\text{$p$}}
\put(14,12){\text{$-$}}
        \end{picture}
    \end{minipage} \!\! \otimes [x]
\right).  
\end{equation}

\appendix
\section{Chain homotopy maps.  }\label{homotopy}
The homotopy connecting $\operatorname{in}$ $\circ$ $\rho_{1'}$ to the identity $h_{1'}$ $:$ $\mathcal{C}\left(~
    \begin{minipage}{30pt}
    \begin{picture}(30,30)
    \put(-2,13.5){$a$}
    \qbezier(6.6,18)(0,25)(0,25)
    \qbezier(0,5)(20,35)(25,18)
    \qbezier(11,14)(20,5)(24.5,13)
    \qbezier(24.5,13)(25.5,15.5)(25,18)
    \end{picture}
    \end{minipage}
\right)$ $\to$ 
$\mathcal{C}\left(~
    \begin{minipage}{30pt}
    \begin{picture}(30,30)
    \put(-2,13.5){$a$}
    \qbezier(6.6,18)(0,25)(0,25)
    \qbezier(0,5)(20,35)(25,18)
    \qbezier(11,14)(20,5)(24.5,13)
    \qbezier(24.5,13)(25.5,15.5)(25,18)
    \end{picture}
    \end{minipage}
\right)$ such that $\delta_{s, t}$ $\circ$ $h_{1'}$ $+$  $h_{1'}$ $\circ$ $\delta_{s, t}$ $=$ $\operatorname{id} - \operatorname{in} \circ \rho_{1'}$, is defined by the formulas: 
\begin{equation}\label{1st}
\begin{split}
\begin{minipage}{30pt}
        \begin{picture}(30,30)
\qbezier(6.6,17)(0,24)(0,24)
\qbezier(0,4)(20,34)(25,17)
\qbezier(11,13)(20,4)(24.5,12)
\qbezier(24.5,12)(25.5,14.5)(25,17)
{\color{red}{\put(12.5,15.5){\circle*{3}}
\put(4.5,15.5){\circle*{3}}
\put(4.5,15.5){\line(1,0){8}}}}
\put(16,14){\text{$p$}}
        \end{picture}
    \end{minipage} \!\! \otimes [xa] &\mapsto~ 
    \begin{minipage}{30pt}
        \begin{picture}(30,30)
\qbezier(6.6,17)(0,24)(0,24)
\qbezier(0,4)(20,34)(25,17)
\qbezier(11,13)(20,4)(24.5,12)
\qbezier(24.5,12)(25.5,14.5)(25,17)
{\color{blue}{\put(8,19.5){\circle*{3}}
\put(8,10.5){\circle*{3}}
\put(8,11.5){\line(0,1){8}}}}
\put(-3.5,12.5){\text{$p$}}
\put(14,12.5){\text{$-$}}
        \end{picture}
    \end{minipage} \!\! \otimes [x],~
    {\text{otherwise}} \mapsto 0.  
\end{split}
\end{equation}
The homotopy connecting ${\operatorname{in}} \circ \rho_{2}$ to the identity $h_{2}$ $:$ $\mathcal{C}\left(~
    \begin{minipage}{30pt}
        \begin{picture}(30,40)
            \qbezier(0,5)(40,20)(0,35)
            \qbezier(11,26)(1,19)(11,14)
            \qbezier(18,29)(24,32)(30,35)
            \qbezier(18,10)(24,7.5)(30,5)
        \end{picture}
    \end{minipage}
~\right) \to \mathcal{C}\left(~
    \begin{minipage}{30pt}
        \begin{picture}(30,40)
            \qbezier(0,5)(40,20)(0,35)
            \qbezier(11,26)(1,19)(11,14)
            \qbezier(18,29)(24,32)(30,35)
            \qbezier(18,10)(24,7.5)(30,5)
        \end{picture}
    \end{minipage}
~\right)$ such that $\delta_{s, t} \circ h_{2}$ $+$ $h_{2} \circ \delta_{s, t}$ $=$ ${\operatorname{id}} - {\operatorname{in}} \circ \rho_{2}$, is defined by the formulas: 
\begin{align}\label{2nd}
&\begin{minipage}{30pt}
        \begin{picture}(30,40)
            \qbezier(0,0)(40,15)(0,30)
            \qbezier(11,21)(1,15)(11,9)
            \qbezier(18,24)(24,27)(30,30)
            \qbezier(18,5)(24,2.5)(30,0)
              {\color{red}{\put(14,28){\circle*{3}}
\put(14,19){\circle*{3}}
\put(14,19){\line(0,1){8}}}}
            {\color{red}{\put(16,7){\circle*{3}}
\put(5,7){\circle*{3}}
\put(5,7){\line(1,0){10}}}}
\put(8,34){$p$}
\put(8,-3){$q$}
        \end{picture}
    \end{minipage} \otimes [xab] \mapsto - \begin{minipage}{30pt}
        \begin{picture}(30,40)
            \qbezier(0,0)(40,15)(0,30)
            \qbezier(11,21)(1,15)(11,9)
            \qbezier(18,24)(24,27)(30,30)
            \qbezier(18,5)(24,2.5)(30,0)
            {\color{blue}
            {\put(8,22.5){\circle*{3}}
\put(20,22.5){\circle*{3}}
\put(9,22.5){\line(1,0){10}}
            }}
            {\color{red}{\put(16,7){\circle*{3}}
\put(5,7){\circle*{3}}
\put(5,7){\line(1,0){10}}}}
\put(5,12){\text{$-$}}
\put(7,32){$p$}
\put(7,-5){$q$}
        \end{picture}
    \end{minipage} \otimes [xb], \quad  \begin{minipage}{30pt}
        \begin{picture}(30,40)
            \qbezier(0,0)(40,15)(0,30)
            \qbezier(11,21)(1,15)(11,9)
            \qbezier(18,24)(24,27)(30,30)
            \qbezier(18,5)(24,2.5)(30,0)
            {\color{blue}
            {\put(8,22.5){\circle*{3}}
\put(20,22.5){\circle*{3}}
\put(9,22.5){\line(1,0){10}}
            }}
            {\color{red}{\put(16,7){\circle*{3}}
\put(5,7){\circle*{3}}
\put(5,7){\line(1,0){10}}}}
\put(5,12){\text{$+$}}
\put(7,32){$p$}
\put(7,-5){$q$}  
        \end{picture}
    \end{minipage}\otimes [xb] \mapsto \begin{minipage}{30pt}
        \begin{picture}(30,40)
            \qbezier(0,0)(40,15)(0,30)
            \qbezier(30,30)(8,20)(30,30)
            \qbezier(11,21)(1,15)(11,9)
            \qbezier(30,0)(10,10)(30,0)
            {\color{blue}
            {\put(8,22.5){\circle*{3}}
\put(20,22.5){\circle*{3}}
\put(9,22.5){\line(1,0){10}}
            }}
            {\color{blue}{\put(10,12){\circle*{3}}
\put(10,3){\circle*{3}}
\put(10,4){\line(0,1){8}}}}
\put(8,30){$p$}
\put(8,-7){$q$}
        \end{picture}
    \end{minipage} \otimes [x], \\ \nonumber
    &\\ \nonumber
     &\text{otherwise} \mapsto 0.  
\end{align}
The homotopy connecting ${\operatorname{in}} \circ \rho_{3}$ to the identity, that is, a map $h_{3} : $ $\mathcal{C}\left(~
    \begin{minipage}{60pt}
\begin{picture}(50,60)
\cbezier(0,0)(0,3)(7,7)(12,12)
\qbezier(19,19)(25,25)(31,31)
\cbezier(38,38)(45,45)(50,50)(50,56)
\cbezier(0,56)(0,50)(17.25,44.25)(34.5,34.5)
\cbezier(34.5,34.5)(44,30)(50,20)(50,0)
\cbezier(17,48)(21,53)(24,54)(24,56)
\cbezier(12,41)(-1,25)(24,10)(24,0)
\put(4,40){$c$}
\put(41,31){$b$}
\put(3,12){$a$}
\end{picture}
\end{minipage}
\right)$ $\to$ $\mathcal{C}\left(~
    \begin{minipage}{60pt}
\begin{picture}(50,60)
\cbezier(0,0)(0,3)(7,7)(12,12)
\qbezier(19,19)(25,25)(31,31)
\cbezier(38,38)(45,45)(50,50)(50,56)
\cbezier(0,56)(0,50)(17.25,44.25)(34.5,34.5)
\cbezier(34.5,34.5)(44,30)(50,20)(50,0)
\cbezier(17,48)(21,53)(24,54)(24,56)
\cbezier(12,41)(-1,25)(24,10)(24,0)
\put(4,40){$c$}
\put(41,31){$b$}
\put(3,12){$a$}
\end{picture}
\end{minipage}\right)$ such that $\delta_{s, t} \circ h_{3}$ $+$ $h_{3} \circ \delta_{s, t}$ $=$ ${\operatorname{id}} - {\operatorname{in}} \circ \rho_{3}$, is defined by the formulas:  
\begin{align}\label{3rd}
&\begin{minipage}{60pt}
\begin{picture}(50,60)
\cbezier(0,0)(0,3)(7,7)(12,12)
\qbezier(19,19)(25,25)(31,31)
\cbezier(38,38)(45,45)(50,50)(50,56)
\cbezier(0,56)(0,50)(17.25,44.25)(34.5,34.5)
\cbezier(34.5,34.5)(44,30)(50,20)(50,0)
\cbezier(17,48)(21,53)(24,54)(24,56)
\cbezier(12,41)(-1,25)(24,10)(24,0)
            {\color{red}{\put(16,21){\circle*{3}}
\put(16,10){\circle*{3}}
\put(16,11){\line(0,1){10}}}}
{\color{blue}
            {\put(4.5,44.5){\circle*{3}}
\put(16.5,44.5){\circle*{3}}
\put(5.5,44.5){\line(1,0){10}}
            }}
            {\color{red}{\put(27,39){\circle*{3}}
\put(27,29){\circle*{3}}
\put(27,29){\line(0,1){10}}}}
            \put(-3,14){$q$}
            \put(45,45){$p$}
            \put(4,48){$r$}
        \end{picture}
    \end{minipage} \otimes [xab] \mapsto -  \begin{minipage}{60pt}
\begin{picture}(50,60)
\cbezier(0,0)(0,3)(7,7)(12,12)
\qbezier(19,19)(25,25)(31,31)
\cbezier(38,38)(45,45)(50,50)(50,56)
\cbezier(0,56)(0,50)(17.25,44.25)(34.5,34.5)
\cbezier(34.5,34.5)(44,30)(50,20)(50,0)
\cbezier(17,48)(21,53)(24,54)(24,56)
\cbezier(12,41)(-1,25)(24,10)(24,0)
            \put(15,28){\text{$-$}}
            {\color{blue}
            {\put(8,44.5){\circle*{3}}
\put(20,44.5){\circle*{3}}
\put(9,44.5){\line(1,0){10}}
            }}
            {\color{blue}{\put(8,15){\circle*{3}}
\put(18,15){\circle*{3}}
\put(9,15){\line(1,0){10}}}}
{\color{red}{\put(27,39){\circle*{3}}
\put(27,29){\circle*{3}}
\put(27,29){\line(0,1){10}}}}
\put(5,-2){$q$}
            \put(45,45){$p$}
\put(5,50){$r$}
        \end{picture}
    \end{minipage} \otimes [xb], 
    \begin{minipage}{60pt}
\begin{picture}(50,60)
\cbezier(0,0)(0,3)(7,7)(12,12)
\qbezier(19,19)(25,25)(31,31)
\cbezier(38,38)(45,45)(50,50)(50,56)
\cbezier(0,56)(0,50)(17.25,44.25)(34.5,34.5)
\cbezier(34.5,34.5)(44,30)(50,20)(50,0)
\cbezier(17,48)(21,53)(24,54)(24,56)
\cbezier(12,41)(-1,25)(24,10)(24,0)
            \put(15,28){\text{$+$}}
            {\color{blue}
            {\put(8,44.5){\circle*{3}}
\put(20,44.5){\circle*{3}}
\put(9,44.5){\line(1,0){10}}
            }}
            {\color{blue}{\put(8,15){\circle*{3}}
\put(18,15){\circle*{3}}
\put(9,15){\line(1,0){10}}}}
{\color{red}{\put(27,39){\circle*{3}}
\put(27,29){\circle*{3}}
\put(27,29){\line(0,1){10}}}}
\put(5,-2){$q$}
\put(40,0){$p$}
\put(5,50){$r$}
        \end{picture}
    \end{minipage}\!\!\!\!\!\! \otimes [xb] \mapsto \begin{minipage}{60pt}
\begin{picture}(50,60)
\cbezier(0,0)(0,3)(7,7)(12,12)
\qbezier(19,19)(25,25)(31,31)
\cbezier(38,38)(45,45)(50,50)(50,56)
\cbezier(0,56)(0,50)(17.25,44.25)(34.5,34.5)
\cbezier(34.5,34.5)(44,30)(50,20)(50,0)
\cbezier(17,48)(21,53)(24,54)(24,56)
\cbezier(12,41)(-1,25)(24,10)(24,0)
            {\color{blue}
            {\put(8,44.5){\circle*{3}}
\put(20,44.5){\circle*{3}}
\put(9,44.5){\line(1,0){10}}
            }}
            {\color{blue}
            {\put(37.5,34.5){\circle*{3}}
\put(25.5,34.5){\circle*{3}}
\put(26.5,34.5){\line(1,0){10}}
            }}
            {\color{blue}{\put(4,15){\circle*{3}}
\put(14,15){\circle*{3}}
\put(5,15){\line(1,0){10}}}}
\put(5,-2){$q$}
\put(40,30){$p$}
\put(5,50){$r$}
        \end{picture}
    \end{minipage} \otimes [x], \\ \nonumber
    \\ \nonumber
&\text{otherwise} \mapsto 0 
\end{align}

\section{Retractions.  }\label{retraction}
The retraction $\rho_{1'} : \mathcal{C}\left(~
    \begin{minipage}{30pt}
        \begin{picture}(30,30)
\qbezier(6.6,18)(0,25)(0,25)
\qbezier(0,5)(20,35)(25,18)
\qbezier(11,14)(20,5)(24.5,13)
\qbezier(24.5,13)(25.5,15.5)(25,18)
        \end{picture}
    \end{minipage}
\right)$ $\to$ 
$ 
\mathcal{C}\left(~\begin{minipage}{30pt}
        \begin{picture}(30,30)
\qbezier(6.6,17)(0,24)(0,24)
\qbezier(0,4)(20,34)(25,17)
\qbezier(11,13)(20,4)(24.5,12)
\qbezier(24.5,12)(25.5,14.5)(25,17)
{\color{blue}{\put(8,19.5){\circle*{3}}
\put(8,10.5){\circle*{3}}
\put(8,11.5){\line(0,1){8}}}}
\put(-3.5,12.5){\text{$p$}}
\put(14,12.5){\text{$+$}}
        \end{picture}
    \end{minipage}\!\! \otimes [x] 
- \operatorname{m}(p:+) \begin{minipage}{30pt}
        \begin{picture}(30,30)
\qbezier(6.6,17)(0,24)(0,24)
\qbezier(0,4)(20,34)(25,17)
\qbezier(11,13)(20,4)(24.5,12)
\qbezier(24.5,12)(25.5,14.5)(25,17)
{\color{blue}{\put(8,19.5){\circle*{3}}
\put(8,10.5){\circle*{3}}
\put(8,11.5){\line(0,1){8}}}}
\put(14,12.5){\text{$-$}}
        \end{picture}
    \end{minipage} \!\! \otimes [x]
    \right)$ is defined by the formulas

\begin{equation}\label{1st-ret}
\begin{split}
&\begin{minipage}{30pt}
        \begin{picture}(30,30)
        \qbezier(6.6,17)(0,24)(0,24)
\qbezier(0,4)(20,34)(25,17)
\qbezier(11,13)(20,4)(24.5,12)
\qbezier(24.5,12)(25.5,14.5)(25,17)
{\color{blue}{\put(8,19.5){\circle*{3}}
\put(8,10.5){\circle*{3}}
\put(8,11.5){\line(0,1){8}}}}
\put(-3.5,12.5){\text{$p$}}
\put(14,12.5){\text{$+$}}
        \end{picture}
    \end{minipage} \!\! \otimes [x] \mapsto \ 
\begin{minipage}{30pt}
        \begin{picture}(30,30)
\qbezier(6.6,17)(0,24)(0,24)
\qbezier(0,4)(20,34)(25,17)
\qbezier(11,13)(20,4)(24.5,12)
\qbezier(24.5,12)(25.5,14.5)(25,17)
{\color{blue}{\put(8,19.5){\circle*{3}}
\put(8,10.5){\circle*{3}}
\put(8,11.5){\line(0,1){8}}}}
\put(-3.5,12.5){\text{$p$}}
\put(14,12.5){\text{$+$}}
        \end{picture}
    \end{minipage}\!\! \otimes [x] 
- \operatorname{m}(p:+) \begin{minipage}{30pt}
        \begin{picture}(30,30)
\qbezier(6.6,17)(0,24)(0,24)
\qbezier(0,4)(20,34)(25,17)
\qbezier(11,13)(20,4)(24.5,12)
\qbezier(24.5,12)(25.5,14.5)(25,17)
{\color{blue}{\put(8,19.5){\circle*{3}}
\put(8,10.5){\circle*{3}}
\put(8,11.5){\line(0,1){8}}}}
\put(14,12.5){\text{$-$}}
        \end{picture}
    \end{minipage} \!\! \otimes [x], \\
&    \begin{minipage}{30pt}
        \begin{picture}(30,30)
\qbezier(6.6,17)(0,24)(0,24)
\qbezier(0,4)(20,34)(25,17)
\qbezier(11,13)(20,4)(24.5,12)
\qbezier(24.5,12)(25.5,14.5)(25,17)
{\color{blue}{\put(8,19.5){\circle*{3}}
\put(8,10.5){\circle*{3}}
\put(8,11.5){\line(0,1){8}}}}
\put(-3.5,12.5){\text{$p$}}
\put(14,12.5){\text{$-$}}
        \end{picture}
    \end{minipage} \!\! \otimes [x],~\   
    \begin{minipage}{30pt}
        \begin{picture}(30,30)
\qbezier(6.6,17)(0,24)(0,24)
\qbezier(0,4)(20,34)(25,17)
\qbezier(11,13)(20,4)(24.5,12)
\qbezier(24.5,12)(25.5,14.5)(25,17)
{\color{red}{\put(12.5,15.5){\circle*{3}}
\put(4.5,15.5){\circle*{3}}
\put(4.5,15.5){\line(1,0){8}}}}
\put(16,14){\text{$p$}}
        \end{picture}
    \end{minipage} \!\! \otimes [xa] \mapsto \ 0.  
\end{split}
\end{equation}
The retraction $\rho_{2} : \mathcal{C}
\left(~
    \begin{minipage}{30pt}
        \begin{picture}(30,30)
            \qbezier(0,0)(40,15)(0,30)
            \qbezier(11,21)(1,15)(11,9)
            \qbezier(18,24)(24,27)(30,30)
            \qbezier(18,5)(24,2.5)(30,0)
            \put(2,20){$a$}
            \put(2,5){$b$}
        \end{picture}
    \end{minipage}
~\right)$ $\to$
$\mathcal{C}\left(~\begin{minipage}{40pt}
        \begin{picture}(30,40)
           \qbezier(0,5)(40,20)(0,35)
            \qbezier(11,26)(1,19)(11,14)
            \qbezier(18,29)(24,32)(30,35)
            \qbezier(18,10)(24,7.5)(30,5)
  {\color{red}{\put(14,33){\circle*{3}}
\put(14,24){\circle*{3}}
\put(14,24){\line(0,1){8}}}}
            {\color{blue}{\put(10,17){\circle*{3}}
\put(10,8){\circle*{3}}
\put(10,9){\line(0,1){8}}}}
\put(-5,19){$p$}
\put(18,19){$q$}
        \end{picture}
    \end{minipage} \otimes [xa] + \begin{minipage}{40pt}
        \begin{picture}(30,40)
            \qbezier(0,5)(40,20)(0,35)
            \qbezier(11,26)(1,19)(11,14)
            \qbezier(18,29)(24,32)(30,35)
            \qbezier(18,10)(24,7.5)(30,5)
            {\color{blue}
            {\put(8,27.5){\circle*{3}}
\put(20,27.5){\circle*{3}}
\put(9,27.5){\line(1,0){10}}
            }}
            {\color{red}{\put(16,12){\circle*{3}}
\put(5,12){\circle*{3}}
\put(5,12){\line(1,0){10}}}}
\put(5,17){\text{$-$}}
\put(1,37){$p:q$}
\put(0,0){$q:p$}
        \end{picture}
    \end{minipage} \otimes [xb]~\right)$ is defined by the formulas 
   \begin{equation}
   \begin{split}
   \begin{minipage}{40pt}
        \begin{picture}(30,40)
            \qbezier(0,5)(40,20)(0,35)
            \qbezier(11,26)(1,19)(11,14)
            \qbezier(18,29)(24,32)(30,35)
            \qbezier(18,10)(24,7.5)(30,5)
  {\color{red}{\put(14,33){\circle*{3}}
\put(14,24){\circle*{3}}
\put(14,24){\line(0,1){8}}}}
            {\color{blue}{\put(10,17){\circle*{3}}
\put(10,8){\circle*{3}}
\put(10,9){\line(0,1){8}}}}
\put(-5,19){$p$}
\put(18,19){$q$}
        \end{picture} 
        \end{minipage} \otimes [xa]
        &\mapsto \begin{minipage}{40pt}
        \begin{picture}(30,40)
            \qbezier(0,5)(40,20)(0,35)
            \qbezier(11,26)(1,20)(11,14)
            \qbezier(18,29)(24,32)(30,35)
            \qbezier(18,10)(24,7.5)(30,5)
  {\color{red}{\put(14,33){\circle*{3}}
\put(14,24){\circle*{3}}
\put(14,24){\line(0,1){8}}}}
            {\color{blue}{\put(10,17){\circle*{3}}
\put(10,8){\circle*{3}}
\put(10,9){\line(0,1){8}}}}
\put(-5,19){$p$}
\put(18,19){$q$}
        \end{picture}
    \end{minipage} \otimes [xa] + \begin{minipage}{40pt}
        \begin{picture}(30,40)
            \qbezier(0,5)(40,20)(0,35)
            \qbezier(11,26)(1,20)(11,14)
            \qbezier(18,29)(24,32)(30,35)
            \qbezier(18,10)(24,7.5)(30,5)
            {\color{blue}
            {\put(8,27.5){\circle*{3}}
\put(20,27.5){\circle*{3}}
\put(9,27.5){\line(1,0){10}}
            }}
            {\color{red}{\put(16,12){\circle*{3}}
\put(5,12){\circle*{3}}
\put(5,12){\line(1,0){10}}}}
\put(5,17){\text{$-$}}
\put(1,37){$p:q$}
\put(0,0){$q:p$}
        \end{picture}
    \end{minipage} \otimes [xb], \\
    \begin{minipage}{40pt}
        \begin{picture}(30,40)
           \qbezier(0,5)(40,20)(0,35)
            \qbezier(11,26)(1,19)(11,14)
            \qbezier(18,29)(24,32)(30,35)
            \qbezier(18,10)(24,7.5)(30,5)
            {\color{blue}
            {\put(8,27.5){\circle*{3}}
\put(20,27.5){\circle*{3}}
\put(9,27.5){\line(1,0){10}}
            }}
            {\color{red}{\put(16,12){\circle*{3}}
\put(5,12){\circle*{3}}
\put(5,12){\line(1,0){10}}}}
\put(5,17){\text{$+$}}
\put(7,37){$p$}
\put(7,0){$q$}
        \end{picture}
    \end{minipage} \otimes [xb] &\mapsto
    - \left(\qquad  \begin{minipage}{40pt}
        \begin{picture}(30,40)
            \qbezier(0,5)(40,20)(0,35)
            \qbezier(11,26)(1,20)(11,14)
            \qbezier(18,29)(24,32)(30,35)
            \qbezier(18,10)(24,7.5)(30,5)
  {\color{red}{\put(14,33){\circle*{3}}
\put(14,24){\circle*{3}}
\put(14,24){\line(0,1){8}}}}
            {\color{blue}{\put(10,17){\circle*{3}}
\put(10,8){\circle*{3}}
\put(10,9){\line(0,1){8}}}}
\put(-23,19){$p:q$}
\put(20,19){$q:p$}
        \end{picture}
    \end{minipage}\quad \otimes [xa] \quad + \quad 
    \begin{minipage}{40pt}
        \begin{picture}(30,40)
            \qbezier(0,5)(40,20)(0,35)
            \qbezier(11,26)(1,19)(11,14)
            \qbezier(18,29)(24,32)(30,35)
            \qbezier(18,10)(24,7.5)(30,5)
            {\color{blue}
            {\put(8,27.5){\circle*{3}}
\put(20,27.5){\circle*{3}}
\put(9,27.5){\line(1,0){10}}
            }}
            {\color{red}{\put(16,12){\circle*{3}}
\put(5,12){\circle*{3}}
\put(5,12){\line(1,0){10}}}}
\put(5,17){\text{$-$}}
\put(-25,42){$(p:q):(q:p)$}
\put(-25,-5){$(q:p):(p:q)$}
        \end{picture}
    \end{minipage} \otimes [xb]~\right), \\
    \text{otherwise} &\mapsto 0.  
\end{split}
\end{equation}
    The retraction $\rho_{3} : \mathcal{C}\left(~\begin{minipage}{60pt}
\begin{picture}(50,60)
\cbezier(0,0)(0,3)(7,7)(12,12)
\qbezier(19,19)(25,25)(31,31)
\cbezier(38,38)(45,45)(50,50)(50,56)
\cbezier(0,56)(0,50)(17.25,44.25)(34.5,34.5)
\cbezier(34.5,34.5)(44,30)(50,20)(50,0)
\cbezier(17,48)(21,53)(24,54)(24,56)
\cbezier(12,41)(-1,25)(24,10)(24,0)
\put(4,40){$c$}
\put(41,31){$b$}
\put(3,12){$a$}
\end{picture}
\end{minipage}\right)$ $\to$
$\mathcal{C}\left(~\begin{minipage}{60pt}
\begin{picture}(50,70)
\cbezier(0,0)(0,3)(7,7)(12,12)
\qbezier(19,19)(25,25)(31,31)
\cbezier(38,38)(45,45)(50,50)(50,56)
\cbezier(0,56)(0,50)(17.25,44.25)(34.5,34.5)
\cbezier(34.5,34.5)(44,30)(50,20)(50,0)
\cbezier(17,48)(21,53)(24,54)(24,56)
\cbezier(12,41)(-1,25)(24,10)(24,0)
            {\color{red}{\put(16,21){\circle*{3}}
\put(16,10){\circle*{3}}
\put(16,11){\line(0,1){10}}}}
{\color{blue}
            {\put(4.5,44.5){\circle*{3}}
\put(16.5,44.5){\circle*{3}}
\put(5.5,44.5){\line(1,0){10}}
            }}
            {\color{blue}
            {\put(32.5,34.5){\circle*{3}}
\put(20.5,34.5){\circle*{3}}
\put(21.5,34.5){\line(1,0){10}}
            }}
            \put(13,14){$q$}
            \put(45,45){$p$}
            \put(4,48){$r$}
        \end{picture}
    \end{minipage} \otimes [xa] + \begin{minipage}{60pt}
\begin{picture}(50,70)
\cbezier(0,0)(0,3)(7,7)(12,12)
\qbezier(19,19)(25,25)(31,31)
\cbezier(38,38)(45,45)(50,50)(50,56)
\cbezier(0,56)(0,50)(17.25,44.25)(34.5,34.5)
\cbezier(34.5,34.5)(44,30)(50,20)(50,0)
\cbezier(17,48)(21,53)(24,54)(24,56)
\cbezier(12,41)(-1,25)(24,10)(24,0)
            \put(15,28){\text{$-$}}
            {\color{blue}
            {\put(8,44.5){\circle*{3}}
\put(20,44.5){\circle*{3}}
\put(9,44.5){\line(1,0){10}}
            }}
            {\color{blue}{\put(8,15){\circle*{3}}
\put(18,15){\circle*{3}}
\put(9,15){\line(1,0){10}}}}
{\color{red}{\put(27,39){\circle*{3}}
\put(27,29){\circle*{3}}
\put(27,29){\line(0,1){10}}}}
\put(-6,-5){$p:q$}
\put(45,0){$q:p$}
\put(5,50){$\tilde{r}$}
        \end{picture}
    \end{minipage}\!\!\!\! \otimes [xb],\right.$ $\left.\begin{minipage}{60pt}
\begin{picture}(50,70)
\cbezier(0,0)(0,3)(7,7)(12,12)
\qbezier(19,19)(25,25)(31,31)
\cbezier(38,38)(45,45)(50,50)(50,56)
\cbezier(0,56)(0,50)(17.25,44.25)(34.5,34.5)
\cbezier(34.5,34.5)(44,30)(50,20)(50,0)
\cbezier(17,48)(21,53)(24,54)(24,56)
\cbezier(12,41)(-1,25)(24,10)(24,0)
            {\color{red}{\put(14,49){\circle*{3}}
\put(14,39){\circle*{3}}
\put(14,39){\line(0,1){10}}}}
        \end{picture}
    \end{minipage} \otimes [x]\right)$ is defined by the formulas 
\begin{align}
&\begin{minipage}{60pt}
\begin{picture}(50,70)
\cbezier(0,0)(0,3)(7,7)(12,12)
\qbezier(19,19)(25,25)(31,31)
\cbezier(38,38)(45,45)(50,50)(50,56)
\cbezier(0,56)(0,50)(17.25,44.25)(34.5,34.5)
\cbezier(34.5,34.5)(44,30)(50,20)(50,0)
\cbezier(17,48)(21,53)(24,54)(24,56)
\cbezier(12,41)(-1,25)(24,10)(24,0)
            {\color{red}{\put(16,21){\circle*{3}}
\put(16,10){\circle*{3}}
\put(16,11){\line(0,1){10}}}}
{\color{blue}
            {\put(4.5,44.5){\circle*{3}}
\put(16.5,44.5){\circle*{3}}
\put(5.5,44.5){\line(1,0){10}}
            }}
            {\color{blue}
            {\put(32.5,34.5){\circle*{3}}
\put(20.5,34.5){\circle*{3}}
\put(21.5,34.5){\line(1,0){10}}
            }}
            \put(13,14){$q$}
            \put(45,45){$p$}
            \put(4,48){$r$}
        \end{picture}
    \end{minipage} \otimes [xa] \mapsto \begin{minipage}{60pt}
\begin{picture}(50,70)
\cbezier(0,0)(0,3)(7,7)(12,12)
\qbezier(19,19)(25,25)(31,31)
\cbezier(38,38)(45,45)(50,50)(50,56)
\cbezier(0,56)(0,50)(17.25,44.25)(34.5,34.5)
\cbezier(34.5,34.5)(44,30)(50,20)(50,0)
\cbezier(17,48)(21,53)(24,54)(24,56)
\cbezier(12,41)(-1,25)(24,10)(24,0)
            {\color{red}{\put(16,21){\circle*{3}}
\put(16,10){\circle*{3}}
\put(16,11){\line(0,1){10}}}}
{\color{blue}
            {\put(4.5,44.5){\circle*{3}}
\put(16.5,44.5){\circle*{3}}
\put(5.5,44.5){\line(1,0){10}}
            }}
            {\color{blue}
            {\put(32.5,34.5){\circle*{3}}
\put(20.5,34.5){\circle*{3}}
\put(21.5,34.5){\line(1,0){10}}
            }}
            \put(13,14){$q$}
            \put(45,45){$p$}
            \put(4,48){$r$}
        \end{picture}
    \end{minipage} \otimes [xa] + \begin{minipage}{60pt}
\begin{picture}(50,70)
\cbezier(0,0)(0,3)(7,7)(12,12)
\qbezier(19,19)(25,25)(31,31)
\cbezier(38,38)(45,45)(50,50)(50,56)
\cbezier(0,56)(0,50)(17.25,44.25)(34.5,34.5)
\cbezier(34.5,34.5)(44,30)(50,20)(50,0)
\cbezier(17,48)(21,53)(24,54)(24,56)
\cbezier(12,41)(-1,25)(24,10)(24,0)
            \put(15,28){\text{$-$}}
            {\color{blue}
            {\put(8,44.5){\circle*{3}}
\put(20,44.5){\circle*{3}}
\put(9,44.5){\line(1,0){10}}
            }}
            {\color{blue}{\put(8,15){\circle*{3}}
\put(18,15){\circle*{3}}
\put(9,15){\line(1,0){10}}}}
{\color{red}{\put(27,39){\circle*{3}}
\put(27,29){\circle*{3}}
\put(27,29){\line(0,1){10}}}}
\put(-6,-5){$p:q$}
\put(45,5){$q:p$}
\put(5,50){$\tilde{r}$}
        \end{picture}
    \end{minipage} \otimes [xb], \\ \nonumber
&\begin{minipage}{60pt}
\begin{picture}(50,70)
\cbezier(0,0)(0,3)(7,7)(12,12)
\qbezier(19,19)(25,25)(31,31)
\cbezier(38,38)(45,45)(50,50)(50,56)
\cbezier(0,56)(0,50)(17.25,44.25)(34.5,34.5)
\cbezier(34.5,34.5)(44,30)(50,20)(50,0)
\cbezier(17,48)(21,53)(24,54)(24,56)
\cbezier(12,41)(-1,25)(24,10)(24,0)
            {\color{red}{\put(14,49){\circle*{3}}
\put(14,39){\circle*{3}}
\put(14,39){\line(0,1){10}}}}
        \end{picture}
    \end{minipage} \otimes [x] \mapsto \begin{minipage}{60pt}
\begin{picture}(50,70)
\cbezier(0,0)(0,3)(7,7)(12,12)
\qbezier(19,19)(25,25)(31,31)
\cbezier(38,38)(45,45)(50,50)(50,56)
\cbezier(0,56)(0,50)(17.25,44.25)(34.5,34.5)
\cbezier(34.5,34.5)(44,30)(50,20)(50,0)
\cbezier(17,48)(21,53)(24,54)(24,56)
\cbezier(12,41)(-1,25)(24,10)(24,0)
            {\color{red}{\put(14,49){\circle*{3}}
\put(14,39){\circle*{3}}
\put(14,39){\line(0,1){10}}}}
        \end{picture}
    \end{minipage} \otimes [x], \\ \nonumber
    \end{align}
    \begin{align*}
&\begin{minipage}{60pt}
\begin{picture}(50,70)
\cbezier(0,0)(0,3)(7,7)(12,12)
\qbezier(19,19)(25,25)(31,31)
\cbezier(38,38)(45,45)(50,50)(50,56)
\cbezier(0,56)(0,50)(17.25,44.25)(34.5,34.5)
\cbezier(34.5,34.5)(44,30)(50,20)(50,0)
\cbezier(17,48)(21,53)(24,54)(24,56)
\cbezier(12,41)(-1,25)(24,10)(24,0)
            \put(15,28){\text{$+$}}
            {\color{blue}
            {\put(8,44.5){\circle*{3}}
\put(20,44.5){\circle*{3}}
\put(9,44.5){\line(1,0){10}}
            }}
            {\color{blue}{\put(8,15){\circle*{3}}
\put(18,15){\circle*{3}}
\put(9,15){\line(1,0){10}}}}
{\color{red}{\put(27,39){\circle*{3}}
\put(27,29){\circle*{3}}
\put(27,29){\line(0,1){10}}}}
\put(5,-2){$q$}
\put(44,0){$p$}
\put(5,50){$r$}
        \end{picture}
    \end{minipage}\!\!\!\!\! \otimes [xb] \mapsto - \begin{minipage}{60pt}
\begin{picture}(50,70)
\cbezier(0,0)(0,3)(7,7)(12,12)
\qbezier(19,19)(25,25)(31,31)
\cbezier(38,38)(45,45)(50,50)(50,56)
\cbezier(0,56)(0,50)(17.25,44.25)(34.5,34.5)
\cbezier(34.5,34.5)(44,30)(50,20)(50,0)
\cbezier(17,48)(21,53)(24,54)(24,56)
\cbezier(12,41)(-1,25)(24,10)(24,0)
            {\color{red}{\put(16,21){\circle*{3}}
\put(16,10){\circle*{3}}
\put(16,11){\line(0,1){10}}}}
{\color{blue}
            {\put(4.5,44.5){\circle*{3}}
\put(16.5,44.5){\circle*{3}}
\put(5.5,44.5){\line(1,0){10}}
            }}
            {\color{blue}
            {\put(32.5,34.5){\circle*{3}}
\put(20.5,34.5){\circle*{3}}
\put(21.5,34.5){\line(1,0){10}}
            }}
            \put(13,14){$p:q$}
            \put(45,45){$q:p$}
            \put(4,48){$r$}
        \end{picture}
    \end{minipage} \otimes [xa] - \begin{minipage}{60pt}
\begin{picture}(50,70)
\cbezier(0,0)(0,3)(7,7)(12,12)
\qbezier(19,19)(25,25)(31,31)
\cbezier(38,38)(45,45)(50,50)(50,56)
\cbezier(0,56)(0,50)(17.25,44.25)(34.5,34.5)
\cbezier(34.5,34.5)(44,30)(50,20)(50,0)
\cbezier(17,48)(21,53)(24,54)(24,56)
\cbezier(12,41)(-1,25)(24,10)(24,0)
            \put(15,28){\text{$-$}}
            {\color{blue}
            {\put(8,44.5){\circle*{3}}
\put(20,44.5){\circle*{3}}
\put(9,44.5){\line(1,0){10}}
            }}
            {\color{blue}{\put(8,15){\circle*{3}}
\put(18,15){\circle*{3}}
\put(9,15){\line(1,0){10}}}}
{\color{red}{\put(27,39){\circle*{3}}
\put(27,29){\circle*{3}}
\put(27,29){\line(0,1){10}}}}
\put(-32,-13){$(q:p):(p:q)$}
\put(30,60){$(p:q):(q:p)$}
\put(5,50){$\tilde{r}$}
        \end{picture}
    \end{minipage} \otimes [xb] \quad - \begin{minipage}{60pt}
\begin{picture}(50,70)
\cbezier(0,0)(0,3)(7,7)(12,12)
\qbezier(19,19)(25,25)(31,31)
\cbezier(38,38)(45,45)(50,50)(50,56)
\cbezier(0,56)(0,50)(17.25,44.25)(34.5,34.5)
\cbezier(34.5,34.5)(44,30)(50,20)(50,0)
\cbezier(17,48)(21,53)(24,54)(24,56)
\cbezier(12,41)(-1,25)(24,10)(24,0)
            {\color{red}
            {\put(15,50){\circle*{3}}
\put(15,40){\circle*{3}}
\put(15,40){\line(0,1){10}}
            }}
            {\color{blue}{\put(8,15){\circle*{3}}
\put(18,15){\circle*{3}}
\put(9,15){\line(1,0){10}}}}
{\color{blue}{\put(32,34){\circle*{3}}
\put(22,34){\circle*{3}}
\put(22,34){\line(1,0){10}}}}
\put(5,-2){$q$}
\put(40,55){$p:r$}
\put(-20,40){$r:p$}
        \end{picture}
    \end{minipage}\!\!\!\! \otimes [xc], \\ \nonumber
&\begin{minipage}{60pt}
\begin{picture}(50,70)
\cbezier(0,0)(0,3)(7,7)(12,12)
\qbezier(19,19)(25,25)(31,31)
\cbezier(38,38)(45,45)(50,50)(50,56)
\cbezier(0,56)(0,50)(17.25,44.25)(34.5,34.5)
\cbezier(34.5,34.5)(44,30)(50,20)(50,0)
\cbezier(17,48)(21,53)(24,54)(24,56)
\cbezier(12,41)(-1,25)(24,10)(24,0)
            {\color{red}{\put(16,21){\circle*{3}}
\put(16,10){\circle*{3}}
\put(16,11){\line(0,1){10}}}}
{\color{blue}
            {\put(4.5,44.5){\circle*{3}}
\put(16.5,44.5){\circle*{3}}
\put(5.5,44.5){\line(1,0){10}}
            }}
            {\color{red}{\put(27,39){\circle*{3}}
\put(27,29){\circle*{3}}
\put(27,29){\line(0,1){10}}}}
            \put(-3,14){$q$}
            \put(45,45){$p$}
            \put(4,48){$r$}
        \end{picture}
    \end{minipage} \otimes [xab] \mapsto \begin{minipage}{60pt}
\begin{picture}(50,70)
\cbezier(0,0)(0,3)(7,7)(12,12)
\qbezier(19,19)(25,25)(31,31)
\cbezier(38,38)(45,45)(50,50)(50,56)
\cbezier(0,56)(0,50)(17.25,44.25)(34.5,34.5)
\cbezier(34.5,34.5)(44,30)(50,20)(50,0)
\cbezier(17,48)(21,53)(24,54)(24,56)
\cbezier(12,41)(-1,25)(24,10)(24,0)
            {\color{red}
            {\put(15,50){\circle*{3}}
\put(15,40){\circle*{3}}
\put(15,40){\line(0,1){10}}
            }}
            {\color{blue}{\put(8,15){\circle*{3}}
\put(18,15){\circle*{3}}
\put(9,15){\line(1,0){10}}}}
{\color{red}{\put(27,39){\circle*{3}}
\put(27,29){\circle*{3}}
\put(27,29){\line(0,1){10}}}}
\put(5,-2){$q$}
\put(45,45){$p$}
            \put(-3,40){$r$}
        \end{picture}
    \end{minipage} \otimes [xbc], \\ \nonumber
\\ \nonumber     
&\text{otherwise} \mapsto 0.  
\end{align*}

\scriptsize{\textsc{Department of Pure and Applied Mathematics Waseda University.  Tokyo 169-8555, Japan.  }}
\scriptsize{{\it E-mail address:} \texttt{noboru@moegi.waseda.jp}}

\begin{thebibliography}{9}
\bibitem{ito3}N. Ito, {\it{On Reidemeister invariance of the Khovanov homology group of the Jones polynomial}}, math.GT/0901.3952.  
\bibitem{ito5}N. Ito, {\it{Chain homotopy maps and a universal differential for Khovanov-type homology}}, math.GT/0907.2104.  
\bibitem{viro} O.Viro, {\it{Khovanov homology, its definitions and ramifications}}, Fund. Math. 184 (2004), 317--342.  
\end{thebibliography}
\end{document}